\documentclass{article}
\begin{document}
\newtheorem{proposition}{Proposition}[section]
\newtheorem{definition}{Definition}[section]
\newtheorem{lemma}{Lemma}[section]
\newcommand{\xl}{\stackrel{\rightharpoonup}{\cdot}}
\newcommand{\xr}{\stackrel{\leftharpoonup}{\cdot}}
\newcommand{\xlplus}{\stackrel{\rightharpoonup}{+}}
\newcommand{\xrplus}{\stackrel{\leftharpoonup}{+}}
\newcommand{\xluplus}{\stackrel{\rightharpoonup}{\uplus}}
\newcommand{\xruplus}{\stackrel{\leftharpoonup}{\uplus}}
\newcommand{\xlodot}{\stackrel{\rightharpoonup}{\odot}}
\newcommand{\xrodot}{\stackrel{\leftharpoonup}{\odot}}
\newcommand{\Ll}{\stackrel{\rightharpoonup}{L}}
\newcommand{\Lr}{\stackrel{\leftharpoonup}{L}}
\newcommand{\Rl}{\stackrel{\rightharpoonup}{R}}
\newcommand{\Rr}{\stackrel{\leftharpoonup}{R}}

\title{\bf Tangent-like Spaces to Local Monoids}
\author{Keqin Liu\\Department of Mathematics\\The University of British Columbia\\Vancouver, BC\\
Canada, V6T 1Z2}
\date{July 19, 2005}
\maketitle

\begin{abstract} \noindent The main new notions are the notions of tangent-like spaces and  local monoids. The main result is the passage from a local monoid to its tangent-like space which is a local Leibniz algebra.\end{abstract}

\bigskip
Based on my belief that Leibniz algebras are too general to establish a fair counterpart of Lie theory in the context of Leibniz algebras, I introduced the notion of a local Leibniz algebra in Section 1.6 of \cite{Liu5}. The purpose of this paper is to construct the analogue of the passage from a linear Lie group to its Lie algebra in the context of local Leibniz algebras. The group-like objects I need in constructing the analogue are local monoids, which are obtained by adding more algebraic structures to monoids with diconjugations introduced in Section 4.2 of \cite{Liu2}. One of the difficulties I experienced in constructing the analogue is to find a suitable definition of a tangent space. The notion of a tangent-like space introduced in this paper is good enough for the purpose of this paper, but some changes may be needed in order to use it to develop the counterpart of differential geometry in a more general context.

\medskip
The paper consists of three sections. Section 1 discusses trimonoids, which give the algebraic foundation of this paper. Section 2 introduces the notion of a local monoid. Section 3 constructs the passage from a local monoid to its tangent-like space which is a local Leibniz algebra.

\medskip
Thoughout this paper, we will use Chapter 3 and Chapter 4 of \cite{Liu2}.

\bigskip
\section{Trimonoids}

We begin this section by introducing the notion of a trisemigroup.

\begin{definition}\label{def1.1} Let $G$ be a nonempty set together with three binary operations: $\,\sharp\,$, $\xl$ and $\xr$. The set $G$ is called a {\bf trisemigroup} if the following two properties hold:
\begin{description}
\item[(i)] For each $\ast\in\{\, \,\sharp\, ,\, \xl  ,\, \xr \, \}$, $(G, \, \ast \,)$ is a semigroup.
\item[(ii)] The three binary operations $\,\sharp\,$, $\xl$ and $\xr$ satisfy the {\bf Hu-Liu triassociative law}:
\begin{equation}\label{eq1} (x\xl y)\,\sharp\, z=x\,\sharp\,(y\xr z)\end{equation}
\begin{equation}\label{eq2} (x\,\sharp\, y)\xl z=x\,\sharp\, (y\xl z)\end{equation}
\begin{equation}\label{eq3} x\xr (y\,\sharp\, z)=(x\xr y)\,\sharp\, z\end{equation}
\begin{equation}\label{eq4} x\xl (y\,\sharp\, z)=x\xl y\xl z\end{equation}
\begin{equation}\label{eq5} (x\,\sharp\, y)\xr z=x\xr y\xr z\end{equation}
for all $x,y,z\in G$.
\end{description}
\end{definition}

A trisemigroup $G$ is also denoted by $(\,G, \, \,\sharp\, , \, \xl, \, \xr \,)$, and the three binary operations $\,\sharp\,$, $\xl$ and $\xr$ are called the {\bf product}, the {\bf left product} and the {\bf right product}, respectively.

\medskip
The Hu-Liu triassociative law was introduced in Definition 3.4 of \cite{Liu2}. If we drop (\ref{eq1}) in Definition~\ref{def1.1}, then we obtain the notion of a {\bf quasitrisemigroup}. The {\bf Hu-Liu quasitriassociative law} means the law consisting of (\ref{eq2}), (\ref{eq3}), (\ref{eq4}) and (\ref{eq5}).

\begin{definition}\label{def1.2} Let $(\,G, \, \,\sharp\, , \, \xl, \, \xr \,)$ be a trisemigroup (quasitrisemigroup). $G$ is called the {\bf trimonoid (quasitrimonoid) with a triunit $e$} if $e\in G$ and
$e$ satisfies
\begin{equation}\label{eq6} 
x\,\sharp\, e=x=e\,\sharp\, x \quad\mbox{for all $x\in G$}.
\end{equation}
and
\begin{equation}\label{eq7} 
x\xl e=x=e\xr x \quad\mbox{for all $x\in G$}.
\end{equation}
\end{definition}

Note that if $e$ is a triunit of a trimonoid $G$, then $e$ satisfies:
\begin{equation}\label{eq8} 
e\xl x=x\xr e \quad\mbox{for all $x\in G$}.
\end{equation}
A triunit of a trimonoid $G$ is also called an {\bf identity} of $G$. By Proposition 3.1 in \cite{Liu2}, the left product $\xl$ and the right product $\xr$ of a trimonoid with a triunit satisfy the diassociative law.

\bigskip
\section{Local Monoids}

This section depends on the notion of a 7-tuple, which was introduced in Definition 4.1 of \cite{Liu2}. Let $(\,A, \,+ , \, \,\sharp\, , \, \xl , \, \xr\,)$ be a 7-tuple with an identity $1^{\times}$. The set 
$$
\hbar ^\times (A):=\{\, e\in A \,|\, \mbox{$e\xr x=x=x\xl e$ for all $x\in A$} \,\}
$$
is called the {\bf halo} of $A$. An element of $\hbar ^\times (A)$ is called a {\bf bar-unit} of $A$. If $e$ is a bar-unit of $A$, then the set
$$
\hbar ^+ (A):=\{\, \alpha \in A \,|\, \mbox{$e\xl \alpha =0$} \,\}
$$
is called the {\bf additive halo} of $A$.

\medskip
By Section 1.2 of \cite{Liu5}, every bar-unit $e$ of a 7-tuple $A$ produces three more binary operations $\xluplus$, $\xruplus$ and $\bullet $ on $A$ in the following way:
\begin{eqnarray}
\label{eq9}x\xluplus y:&=&x+e\xl y,\\
\label{eq10}x\xruplus y:&=&x\xr e+ y,\\
\label{eq11} x\bullet  y: &=& x\xr y +x\xl y - x\xr e\xl y,
\end{eqnarray}
where $x$, $y\in  R$. The binary operations $\xluplus$ and $\xruplus$ are called the {\bf left addition} and the {\bf right addition} induced by $e$, respectively. The binary operation $\bullet $ is called the {\bf Hu-Liu product} induced by $e$. Hence, a 7-tuple 
$(\,A, \,+ , \, \,\sharp\, , \, \xl , \, \xr\,)$ always carries the following seven binary operations:
\begin{equation}\label{eq12}
+ , \quad  \,\sharp\, , \quad \xl , \quad \xr , \quad \xluplus , \quad \xruplus , \quad \bullet .
\end{equation}
This is our reason of using 7-tuple to name the algebraic object introduced in Definition 4.1 of \cite{Liu2}. If it is necessary to indicate explicitely that the left addition, the right addition and the Hu-Liu product are induced by a bar-unit $e$, then we use $\xluplus _e$,
$\xruplus _e$ and $\bullet _e$ to denote $\xluplus$, $\xruplus$ and $\bullet$, respectively.

\begin{proposition}\label{pr2.1} Let $(\,A, \,+ , \, \,\sharp\, , \, \xl , \, \xr\,)$ be a 
7-tuple. If $e$ is a bar-unit of $A$, then $(\,A, \, \bullet _e , \, \xl , \, \xr\,)$ is a quasitrimonoid with the triunit $e$.
\end{proposition}

\medskip
\noindent
{\bf Proof} By Proposition 1.3 in \cite{Liu5}, $(\,A, \, \bullet _e \,)$ is a monoid with the unit $e$. A direct computation shows that the three binary operations $\bullet _e$, $\xl$ and $\xr$ satisfy the Hu-Liu quasitriassociative law. Hence, Proposition~\ref{pr2.1} is true.

\hfill\raisebox{1mm}{\framebox[2mm]{}}

\bigskip
The next definition is based on the equation (6.1) in \cite {Liu2}.

\begin{definition}\label{def2.1}  An element $x$ of a 7-tuple 
$(\,A, \,+ , \, \,\sharp\, , \, \xl , \, \xr\,)$ is said to be {\bf one-sided invertible} if there exist two elements $x_{e}^{\stackrel{\ell}{-}1}$ and $x_{e}^{\stackrel{r}{-}1}$ of $A$ such that
\begin{equation}\label{eq13} 
x_e^{\stackrel{\ell}{-}1}\xl x=e=x\xr x_e^{\stackrel{r}{-}1} 
\quad\mbox{for some $e\in \hbar ^\times (A)$}.
\end{equation}
\end{definition}

The two elements 
$x_e^{\stackrel{\ell}{-}1}$ and $x_e^{\stackrel{r}{-}1}$ are called the {\bf left inverse} and the {\bf right inverse} of $x$ with respect to the bar-unit $e$, respectively. By Proposition 6.2 in \cite{Liu2},  $x_e^{\stackrel{\ell}{-}1}=x_e^{\stackrel{r}{-}1}$ if and only if $e$ satisfies  (\ref{eq8}). The set of all 
one-sided invertible elements of a 7-tuple $A$ is denoted by 
$$ A^{\check{-}1}: =\{\, a\in A \,|\,\,\mbox{$y\xl a=e=a\xr z$ for some $y, z\in A$} \,\},$$
where $e$ is a fixed bar-unit of $A$. The definition of $ A^{\check{-}1}$ does not depend on the choice of the bar-unit $e$.

\medskip
For an one-sided invertible element $a$ of a 7-tuple $A$, we define a map
$\Psi _a : A\to A$ by
\begin{equation}\label{eq14}
\Psi _a (x): =a_{e}^{\stackrel{\ell}{-}1}\xr x\xl a=
a_{e}^{\stackrel{r}{-}1}\xr x\xl a\quad\mbox{for $x\in A$},
\end{equation}
where $e$ is a bar-unit of $A$. The map $\Psi _a$ is called the {\bf diconjugation} of $A$ determined by $a$. Since the definition of $\Psi _a $ is independent of the choice of the bar-unit $e$, (\ref{eq14}) is also written as
$$
\Psi _a (x): =a^{\stackrel{\ell}{-}1}\xr x\xl a=
a^{\stackrel{r}{-}1}\xr x\xl a\quad\mbox{for $x\in A$},
$$
where $a^{\stackrel{\ell}{-}1}$ and $a^{\stackrel{r}{-}1}$ denote the left inverse and the right  inverse of $a$ with respect to any bar-unit of $A$, respectively.

\begin{proposition}\label{pr2.2} Let $(\,A, \,+ , \, \,\sharp\, , \, \xl , \, \xr\,)$ be a 7-tuple with an identity $1^\times$.
\begin{description}
\item[(i)] $(\,\hbar ^\times (A), \, \,\sharp\, \,)$ is a monoid with the unit $1^\times$.
\item[(ii)] $(\,\hbar ^+ (A), \, + , \, \,\sharp\, \,)$ is a rng.
\item[(iii)] For any $e\in \hbar ^\times (A)$, $(\,\hbar ^\times (A), \, \bullet _e \,)$ is a monoid with the unit $e$.
\item[(iv)] For any $a\in A^{\check{-}1}$, we have
$$
\Psi _a \Big( \hbar ^\times (A) \Big)=\hbar ^\times (A) \quad\mbox{and}\quad 
\Psi _a \Big( \hbar ^+ (A) \Big)=\hbar ^+ (A).
$$
\item[(v)] If $a\in A^{\check{-}1}$, then the diconjugation $\Psi _a$ preserves each of the first four binary operations in the list (\ref{eq12}) and
\begin{eqnarray*}
\Psi _a \Big( x \xluplus _e y)&=&\Psi _a (x) \xluplus _{\Psi _a (e)} \Psi _a (y),\\ 
\Psi _a \Big( x \xruplus _e y)&=&\Psi _a (x) \xruplus _{\Psi _a (e)} \Psi _a (y),\\ 
\Psi _a \Big( x \bullet _e y)&=&\Psi _a (x) \bullet _{\Psi _a (e)} \Psi _a (y), 
\end{eqnarray*}
where $x$, $y\in A$ and $e\in \hbar ^\times (A)$.
\end{description}
\end{proposition}

\medskip
\noindent
{\bf Proof} They are the direct consequences of the Hu-Liu triassociative law.

\hfill\raisebox{1mm}{\framebox[2mm]{}}

\bigskip
An element $1^{\sharp}$ of a 7-tuple $(\,A, \,+ , \, \,\sharp\, , \, \xl , \, \xr\,)$ is called the {\bf local identity} if $(\,\hbar ^+ (A), \, + , \, \,\sharp\, \,)$ is a ring with the identity $1^{\sharp}$. The notion of a local identity is of importance to rewriting commutative ring theory in a more general context (see Chapter 4 or Chapter 5 in \cite{Liu5} for the application of local identity).

\medskip
Let $(\,A, \,+ , \, \,\sharp\, , \, \xl , \, \xr\,)$ be a 7-tuple with an identity $1^\times$. An element $a$ of $A$ is said to be {\bf invertible} if there exists an element $b$ of $A$ such that
\begin{equation}\label{eq15}
a\,\sharp\, b=1^\times = b\,\sharp\, a.
\end{equation}
The element $b$ satisfying (\ref{eq15}) is called the {\bf inverse} of $a$ and is denoted by $a^{-1}$. We use $A^{-1}$ and $\hbar ^\times (A)^{-1}$ to denote the set of all invertible elements in $A$ and $\hbar ^\times (A)$, respectively. 

\medskip
We now define local monoids, which are the group-like object we need to construct the analogue of the passage from a Lie group to its Lie algebra in the context of local Leibniz algebras. 

\begin{definition}\label{def2.2} Let $(\,A, \,+ , \, \,\sharp\, , \, \xl , \, \xr\,)$ be a 
7-tuple with an identity $1^\times$. A subset $G$ of $A^{\check{-}1}$ is called a 
$(\,\Delta , \, \Omega \, )$-{\bf local monoid} of $A$ if the following five properties hold:
\begin{description}
\item[(i)] $\Psi _a (G)\subseteq G$ for each $a\in G$.
\item[(ii)] $\Delta\subseteq \{\, \,\sharp\, ,\, \bullet_e \,|\, \mbox{$e\in 
G\cap \hbar ^\times (A)$ and $\bullet _e$ is the Hu-Liu product induced by $e$}\,\}$.
\item[(iii)] $\bullet _e\in \Delta$ $\Longrightarrow$ $\bullet _{\Psi _a (e)}\in \Delta$ for each
$a\in G$.
\item[(iv)] $(\,G, \, \ast \, \,)$ is a monoid with a unit for each $\ast \in \Delta$.
\item[(v)] $\Omega\subseteq G\cap \hbar ^\times (A)^{-1}$, $\Psi _a (\Omega)\subseteq \Omega$ for each $a\in G$ and $(\,\Omega , \, \,\sharp\, \,)$ is a group with the identity $1^\times$.
\end{description}
\end{definition}

It is clear that if $G$ is a $(\,\Delta , \, \Omega \, )$-local monoid of a 7-tuple $A$, then $G$ is also a $(\,\Delta , \, \{1^\times\} \, )$-local monoid of $A$.
We will see in the next section that every $(\,\Delta , \, \Omega \, )$-local monoid $G$ of a finite dimensional complete 7-tuple $A$ produces a local Leibniz algebra, where the notion of a finite dimensional complete 7-tuple was introduced in Definition 4.4 of \cite{Liu2}.

\bigskip
\section{Tangent-like Spaces}

The notion of a local Leibniz algebra was introduced in Definition 1.13 of \cite{Liu5} by using the Leibniz identity, the Jacobi identity and the Hu-Liu identity. Before presenting the definition of a local Leibniz algebra, we explain where the Hu-Liu identity comes from.

\medskip
In a 7-tuple $(\,A, \,+ , \, \,\sharp\, , \, \xl , \, \xr\,)$, we can define a angle bracket $\langle \, , \, \rangle $ and a square bracket $[ \, , \, ]$ by
\begin{equation}\label{eq16} 
\langle x , y\rangle := x\xl y - y\xr x \quad\mbox{for $x$, $y\in A$}
\end{equation}
and
\begin{equation}\label{eq17} 
[x ,\, y] := x \,\sharp\, y - y \,\sharp\, x
\quad\mbox{ for $x$, $y\in A$.}
\end{equation}
By Proposition 4.8 in \cite{Liu2}, the angle bracket $\langle \, ,\, \rangle$ defined by (\ref{eq16}) and the square bracket $[\,, \,]$ defined by (\ref{eq17}) satisfy the following identity:
\begin{eqnarray}\label{eq18}
&&[x, \langle y, z \rangle ] +[y, \langle z, x \rangle ] +
[z, \langle x, y \rangle ] \nonumber\\
&=&\;\; [x, \langle z, y \rangle ] +[z, \langle y, x \rangle ]+[y, \langle x, z \rangle ] +\nonumber\\
 &&+\langle [ x, y], z \rangle + \langle [ y, z], x \rangle + \langle [ z, x], y \rangle ,
\end{eqnarray}
where $x$, $y$, $z\in A$. If $x:=\alpha$ and $y:=\beta$ are two elements of the additive halo
$\hbar ^+(A)$, then (\ref{eq18}) becomes
\begin{equation}\label{eq19} 
\langle[\alpha , \beta ] , z \rangle+[\langle \beta , z \rangle , \alpha ]+
[\beta , \langle \alpha , z \rangle ]=0
\end{equation}
for all $\alpha$, $\beta \in \hbar ^+(A)$ and $z\in A$. The identity (\ref{eq19}) is called the {\bf Hu-Liu identity}.

\begin{definition}\label{def3.1} A vector space $L$ over a field ${\bf k}$ is called a {\bf local Leibniz Algebra} if there exists a subspace $L_1$ of $L$ such that the following three properties hold:
\begin{description}
\item[(i)] There is a bilinear map $\langle \, , \, \rangle : L\times L \to L$ satisfying the Leibniz identity and
\begin{equation}\label{eq20} 
\langle L , \, L_1 \rangle =0, \quad 
\langle L_1 , \, L \rangle \subseteq L_1 .
\end{equation}
\item[(ii)] There is a bilinear map $[ \, , \, ] : L_1\times L_1 \to L_1$ satisfying the Jacobi identity and
\begin{equation}\label{eq21} 
[\alpha , \beta]= -[\beta , \alpha ],
\end{equation}
where $\alpha$, $\beta$, $\gamma\in L_1$.
\item[(iii)] The angle bracket $\langle \, , \, \rangle $ and the square bracket $[ \, , \, ]$ satisfy the Hu-Liu identity (\ref{eq19})
for all $\alpha$, $\beta \in L_1$ and $z\in L$.
\end{description}
\end{definition}

The subspace $L_1$ in Definition~\ref{def3.1} is called the {\bf local part} of the local Leibniz algebra $L$.

\medskip
Clearly, a 7-tuple $(\,A, \,+ , \, \,\sharp\, , \, \xl , \, \xr\,)$ is a trialgebra introduced by Definition 1.12 in \cite{Liu5}. By Proposition 1.11 in  \cite{Liu5}, 
$(\, A, \, + , \langle \, , \, \rangle , \, [ \, , \, ] \,)$ is a local Leibniz Algebra with the local part $\hbar ^+(A)$, where the angle bracket $\langle \, , \, \rangle $ and the square bracket $[ \, , \, ]$ are defined by (\ref{eq16}) and (\ref{eq17}), respectively. This fact gives the passage from a 7-tuple to a local Leibniz Algebra. 

\medskip
In the remaining part of this paper, ${\bf k}$ will denote the field $\mathcal{R}$ of real numbers or the field $\mathcal{C}$ of complex numbers.

\begin{definition}\label{def3.2} Let $(\,A, \,+ , \, \,\sharp\, , \, \xl , \, \xr\,)$ be a finite dimensional complete 7-tuple with an identity $1^\times$. If $G$ is a 
$(\,\Delta , \, \Omega \, )$-local monoid of $A$. The {\bf tangent-like space} 
$T_{\Delta, \Omega} (G)$ to $G$ is defined by
\begin{equation}\label{eq22} 
T_{\Delta, \Omega} (G):= T_\Omega (G)+ \sum _{\ast\in \Delta}T_\ast (G),
\end{equation}
where
$$
T_\Omega (G):=
\{\, A'(0) \, | \,\mbox{$A(t)$ is a differentiable curve in $\Omega$ with $A(0)=1^\times$ }\},
$$
$$
T_{\bullet _e}(G):=\{\, a'(0) \,|\, \mbox{$a(t)$ is a differentiable curve in $G$ with $a(0)=e$} \, \}
$$
and
$$
T_{\,\sharp\,}(G):=\{\, a'(0) \,|\, \mbox{$a(t)$ is a differentiable curve in $G$ with $a(0)=1^\times$} \, \}.
$$

\end{definition}

The next proposition gives the main result of this paper.

\begin{proposition}\label{pr3.1} Let $(\,A, \,+ , \, \,\sharp\, , \, \xl , \, \xr\,)$ be a finite dimensional complete 7-tuple with an identity $1^\times$. If $G$ is a 
$(\,\Delta , \, \Omega \, )$-local monoid of $A$, then the tangent-like space 
$(\, T_{\Delta, \Omega} (G), \, + , \langle \, , \, \rangle , \, [ \, , \, ] \,)$ to $G$ is a real local Leibniz algebra with the local part $T_\Omega (G)$, where the angle bracket $\langle \, , \, \rangle $ and the square bracket $[ \, , \, ]$ are defined by (\ref{eq16}) and (\ref{eq17}), respectively. Moreover, the local part $T_\Omega (G)$ is a subset of the additive halo $\hbar ^+(A)$ of $A$.
\end{proposition}

\medskip
\noindent
{\bf Proof} First, using the standard argument in linear Lie groups, we have
\begin{equation}\label{eq23} 
\mbox{$(\, T_\Omega (G), \, + , \, [ \, , \, ] \,)$ is a real Lie algebra}
\end{equation}
and
\begin{equation}\label{eq24} 
\mbox{$(\, T_\ast (G), \, + , \,)$ is a real vector space for each 
$\ast\in \Delta$.}
\end{equation}

\medskip
Next, using the argument in the proof of Proposition 4.4 of \cite{Liu2}, we have
\begin{equation}\label{eq25} 
\langle \, T_\Omega (G) , \, T_\Omega (G)\, \rangle =0,
\end{equation}
\begin{equation}\label{eq26} 
\langle \, T_\ast(G) , \, T_\Omega (G)\, \rangle =0 
\quad\mbox{for $\ast\in \Delta$},
\end{equation}
\begin{equation}\label{eq27} 
\langle \, T_\Omega (G) , \, T_\ast(G)\, \rangle \subseteq T_\Omega (G)
\quad\mbox{for $\ast\in \Delta$},
\end{equation}
\begin{equation}\label{eq28} 
\langle \, T_{\ast_1}(G) , \, T_{\ast_2}(G)\, \rangle \subseteq 
\sum _{\ast\in \Delta}T_\ast (G)
\quad\mbox{for $\ast_1$, $\ast_2\in \Delta$}.
\end{equation}

By (\ref{eq23}) and (\ref{eq24}), $T_{\Delta, \Omega} (G)$ is a real vector space, 
$T_\Omega (G)$ is a subspace of $T_{\Delta, \Omega} (G)$, and 
$T_\Omega (G)$ is closed under the square bracket $[ \, , \, ]$. By (\ref{eq25}), (\ref{eq26}), (\ref{eq27}) and (\ref{eq28}), $T_{\Delta, \Omega} (G)$ is closed under the angle bracket $\langle \, , \, \rangle $ and (\ref{eq20}) holds for 
$L_1:=T_\Omega (G)$. This proves that $T_{\Delta, \Omega} (G)$ is a real local Leibniz algebra with the local part $T_\Omega (G)$.

\medskip
Finally, since $T_\Omega (G)$ is a closed real subspace of $A$ and 
$$\hbar ^\times (A)-\hbar ^\times (A)\subseteq \hbar ^+ (A),$$ 
we have that $T_\Omega (G)$ is a real subspace of $\hbar ^+ (A)$.

\hfill\raisebox{1mm}{\framebox[2mm]{}}

\bigskip
Note that if $G$ is a $(\,\Delta , \, \{1^\times\}\, )$-local monoid, then the tangent-like space $T_{\Delta, \{1^{\times}\}}(G)$ to $G$ is a Leibniz algebra. Moreover, if $(\,G, \, \,\sharp\, \,)$ is a monoid with a unit $1^\times$, then $G$ is a  
$(\,\{\,\sharp\,\} , \, \{1^\times\} \, )$-local monoid and the tangent-like space 
$T_{\{\,\sharp\,\}, \{1^{\times}\}}(G)=T_{\,\sharp\,}(G)$ is a Leibniz algebra, which is Proposition 4.4 of \cite{Liu2}. Hence, the passage established in Proposition 4.4 of \cite{Liu2} is contained in Proposition~\ref{pr3.1}. 

\medskip
The main ideal of this paper can be used to construct the analogue of the passage from a linear Lie group to its Lie algebra in the context of Hu-Liu Leibniz algebras, where Hu-Liu Leibniz algebras are more general than local Leibniz algebras and were introduced in Definition 1.15 of \cite{Liu5}.


\bigskip
\begin{thebibliography}{99}
\bibitem{Kinyon} Michael K. Kinyon, \textsl{Leibniz Algebras, Lie Racks, and digroups}, arXiv: math. RA/0403509v2 31 Mar 2004.
\bibitem{Liu2} Keqin Liu, \textsl{The Generalizations of Groups}, \quad Research Monographs in Mathematics {\bf 1}, 153 Publishing (http://web.ncf.ca/ee722/153publishing/), 2004.
\bibitem{Liu5} Keqin Liu, \quad\textsl{Introduction to Trirings}, \quad Research Monographs in Mathematics {\bf 2}, 153 Publishing (http://web.ncf.ca/ee722/153publishing/), 2005.
\end{thebibliography}
\end{document}